\documentclass{statsoc}
\usepackage{graphicx}
\usepackage{epsfig}
\usepackage[usenames]{color}
\usepackage{epsfig, bm}
\usepackage{latexsym}
\usepackage{bbding}
\usepackage{amssymb}
\usepackage{ifsym}
\usepackage{wasysym}

\RequirePackage{color}
\definecolor{MyDarkGreen}{rgb}{0.02,0.60,0.06}

\title[Critical mass of research groups]{Statistics of statisticians: Critical mass of statistics and operational research groups in the UK}
\author[Kenna and Berche]{Ralph Kenna} 
\address{Applied Mathematics Research Centre,
Coventry University,
Coventry, CV1 5FB, England.}
\email{r.kenna@coventry.ac.uk}
\author[Kenna and Berche]{Bertrand Berche}
\address{Statistical Physics Group,
 Institut Jean Lamour\footnote{Laboratoire associ\'e au CNRS UMR 7198},
 CNRS -- Nancy Universit\'{e} -- UPVM, B.P. 70239,
 F -- 54506 Vand{\oe}uvre l\`es Nancy Cedex, France}
\begin{document}
\maketitle

%-----------------------------------------------------------------------
{\Large
  \begin{abstract}
%-----------------------------------------------------------------------
Using a recently developed model, inspired by mean field theory in statistical physics, and data from the UK's Research Assessment Exercise,
we analyse the relationship between the quality of statistics and operational research groups and the quantity researchers in them.
Similar  to other academic disciplines, we provide evidence for a linear dependency of quality on quantity up to an upper critical mass,
which is interpreted as the average maximum number of colleagues with whom a researcher can communicate meaningfully within a research group.
The model also predicts a lower critical mass, which research groups should strive to achieve to avoid extinction.
For statistics and operational research, the lower critical mass is estimated to be $9 \pm 3$. 
The upper critical mass, beyond which research quality does not significantly depend on group size, is about twice this value. 
%-----------------------------------------------------------------------
  \end{abstract} }
%-----------------------------------------------------------------------
%
%  \thispagestyle{empty}
%
%***********************************************************************
%
%  \newpage
%
%-----------------------------------------------------------------------
%                  \pagenumbering{arabic}
%-----------------------------------------------------------------------

%%%%%%%%%%%%%%%%%%%%%%%%%%%%%%%%%%%%%%%%%%%%%%%%%%%%%%%%%%%%%%%%%%%
\section{Introduction}
%%%%%%%%%%%%%%%%%%%%%%%%%%%%%%%%%%%%%%%%%%%%%%%%%%%%%%%%%%%%%%%%%%%

The notion of {\emph{critical mass}\/} in research has been around for a long time without proper definition. 
As governments, funding councils and universities seek indicators to measure research quality and to pursue greater 
efficiencies in the research sector, critical mass is becoming an increasingly important concept at managerial and policy-making level.
However, until very recently there have been no successful attempts to quantify this notion (Harrison, 2009).
It has been described by Evidence (2010) as ``some minimum size threshold for effective performance''  and, as such, has been 
linked to the idea that benefit accrues through increase of scale of research groups. However, although Evidence (2010) 
demonstrated ``a relationship of some kind between larger 
units and relatively high citation impact'', indications of such a threshold have been lacking.

We recently presented a model for the relationship between quality of research groups and their quantity (Kenna and Berche, 2010a). 
This model was inspired by mean-field theories of statistical physics and 
allowed for a quantitative definition of critical mass. In fact there are two critical masses in research and their values are discipline dependent. 
Instead of a threshold group size  above which research quality improves, 
we have shown that there is a breakpoint or {\emph{upper critical mass}\/} beyond which the linear dependency of research quality on group quantity reduces.
%This breakpoint value,  which we denote by $N_c$, is dependent on the academic discipline under consideration.
Denoting this  value by $N_c$, we showed that the strength of the overall research sector in a given discipline is improved
by supporting groups whose size are less than $N_c$, provided they are  bigger than a second critical mass, which we denote by $N_k$.
Groups whose size are smaller than $N_k$ are vunerable and should seek to achieve the lower critical mass for long-term viability. 
The two critical masses are related by a scaling relation,
\begin{equation}
 N_c = 2 N_k\,.
 \label{scaling}
\end{equation}
We classify research groups of size $N$ within a given discipline as small, medium and large according to whether $N < N_k$, $N_k \le N < N_c$ or $N \ge Nc$, respectively.

We recently determined the critical masses of a multitude of academic disciplines by applying statistical analyses to the results
of the UK's most recent {\emph{Research Assessment Exercise}\/} (RAE) in which the quality of research groups were measured (Kenna and Berche, 2010b).
Notably absent from our analaysis, however, were the statistics and operational research groups, as these were less straightforward to analyse than 
other subject areas. 
Here we rectify this omission  by a careful analysis of these disciplines. 
Our main result is that the lower critical mass, which statistics and operational research groups should attain to be viable in the long term, is 
\begin{equation}
 N_k = 9 \pm 3\,.
 \label{result}
 \end{equation}
 
 In Section~2 we summarize our model and how we derive critical masses from it. We also discuss the research assessment exercise. In Section~3 we apply the model and  statistical analysis  to the results of the RAE for statistics and operational research groups.  We conclude in Section~4, where implications for policy and management are briefly  discussed.

%%%%%%%%%%%%%%%%%%%%%%%%%%%%%%%%%%%%%%%%%%%%%%%%%%%%%%%%%%%%%%%%%%%
\section{Quality and quantity in research}
%%%%%%%%%%%%%%%%%%%%%%%%%%%%%%%%%%%%%%%%%%%%%%%%%%%%%%%%%%%%%%%%%%%

Our model is based on the idea that research groups are {\emph{complex systems}\/}, for which the properties 
of the whole are not simple sums of the corresponding properties of the individual parts.
Instead, interactions between individuals within research groups have to also be taken into account. 
The strength of an individual within a research group is a function of many factors: their intrinsic calibre and training,
their teaching and administrative loads, library facilities, journal access, extramural collaboration, the quality of management, 
and even confidence gained by previous successes as well as the prestige of the institution and other factors. 
We denote the average individual research strength within the $g^{\rm{th}}$ research group in a given academic discipline, 
resulting from all of these (and any other) factors  by $a$. 
The overall calibre of a research group comprising $N$ individuals  is also dependent on the extent of, and strength of, the communication links  between them.
We denote the average strength of the $N(N-1)/2$ interactions between the $N$ individuals in the $g^{\rm{th}}$ group by $b$.
The overall strength of the group is therefore given by 
\begin{equation}
 S = Na + \frac{1}{2} N(N-1) b \,.
 \label{Sg1}
\end{equation}

However, once the size of a research group becomes too large (say above a cutoff value $N_c$), meaningful communication between {\emph{all}\/}
pairs of individuals becomes impossible. In this case, the group may fragment into ${\cal{N}}$  subgroups, of average size $M = N/{\cal{N}}$, say. 
If the average strength of interaction between the  subgroups is $c$,  the overall strength of the group becomes
\begin{equation}
 S = Na + \frac{1}{2} N(M -1) b + \frac{1}{2} {\cal{N}}({\cal{N}}-1) c\,.
 \label{Sg2}
\end{equation}

We denote by $\langle{S}\rangle$ the {\emph{expected}\/} strength of a group of size $N$ and we define the quality of such a research group to be the average strength per head:
\begin{equation}
 s = \frac{S}{N}\,.
 \label{quality}
\end{equation}
Gathering terms of the same order in $N$, we arrive at a form for the {\emph{expected}\/} dependency of research-group quality on research-group quantity,
\begin{equation}
 \langle{s}\rangle = \left\{ \begin{array}{ll}
             a_1 + b_1 N &  {\mbox{if $N \le N_c$}} \\
             a_2 + b_2 N &  {\mbox{if $N \ge N_c$}}.
             \end{array}
     \right.
\label{Nc}
\end{equation}
%where the $a_i$ and $b_i$ are related to the parameters appearing in (\ref{Sg1}) and (\ref{Sg2}) for $i=1,2$, respectively.

We considered the effect on the overall strength of a discipline by adding new researchers (Kenna and Berche, 2010a). Asking the question whether
it is better, on average, to allocate new researchers to a group with $N>N_c$ or $N<N_c$ members, we  found that the latter is preferable
provided $N>N_k$, where $N_k$ is given by Eq.(\ref{scaling}). This is equivalent to maximising the gradient of the strength function $\langle{S (N)}\rangle$.
We also considered the consequences of transferring researchers from large to small/medium groups and found that such a movement is expected to be 
beneficial to society as a whole, provided the recipient group is not too small (i.e., provided, again, that it has over $N_k$ members).
Thus there are two critical masses in research, which we name {\emph{lower}\/} ($N_k$)  and {\emph{upper}\/} ($N_c$). Of these, the former corresponds more closely to the
traditional, intuitive notion of critical mass, although there is no threshold value beyond which research quality suddenly improves (Evidence, 2010).

To implement the model (\ref{Nc}), we require a set of empirical data on the quality and quantity of research groups. 
The RAE is an evaluation process undertaken approximately every 5 years on behalf of the 
funding bodies for universities in the UK. The results of the RAE are used to allocate funding to such  higher education institutes for the subsequent years. 
The last RAE was carried out in 2008.
Research groups were examined to  determine the proportion of research submitted categorized as follows:
\begin{itemize}
\item
4*: Quality that is world-leading in terms of originality, significance and rigour 
\item
3*: Quality that is internationally excellent in terms of originality, significance and rigour but which nonetheless falls short of the highest standards of excellence 
\item
2*: Quality that is recognised internationally in terms of originality, significance and rigour 
\item
1*: Quality that is recognised nationally in terms of originality, significance and rigour
\item Unclassified: Quality that falls below the standard of nationally
 recognised work.
\end{itemize}
A formula is then used to determine how funding is distributed to research groups. 
The 2009 formula used by the Higher Education Funding Council for England  weighs each
rank in such a way that 4* and 3* research  respectively receive seven and three times the amount of funding 
allocated to 2* research, and 1* and unclassified research attract no funding.
This funding formula may therefore be considered to represent a measurement of quality of each research group.
(In 2010, after lobbying by the larger, research intensive universities  the English funding formula was 
changed so that 4*  research  receives nine times the funding  allocated to 2* research. 
We have  checked that the 2010 formula produces no significant 
change to the results presented here.)

%The higher-education community of universities in the UK is organised into a number of represenation or mission groups.
%The {\emph{Russell Group}\/} and {\emph{1994 Group}\/}  consist of research intensive  universities mostly  with and without medical schools, respectively. 
%The {\emph{Million+ Group}\/} of modern universities mostly represents those founded after 1992 while the 
%{\emph{University Alliance}} is of business facing universities. 
%None of the {\emph{Guild~HE}} universities, whose focus is teaching, submitted to the statistics Unit of Assessment and the remaining universities are 
%unaffiliated. [SHOULD WE MENTION THE MISSION GROUPS AT ALL OR IS THIS JUST DILUTING THE MESSAGE?]

From the outset, we acknowledge that there are obvious assumptions underlying our analysis and limits to what can be achieved. 
Firstly, we use the term ``group'' in the sense of RAE. This means the collection of staff included in a submission to one of the 67 {\emph{Units of Assessment}\/} (UOA's). 
RAE groups are not always identical to  administrative departments within universities, but we assume that they represent a coherent group for research purposes.
Individuals submitted to RAE are drawn from academic staff who were in post and on the payroll of the submitting higher education
institution on the census date (31~October~2007). We assume that the RAE process is fair and unbiased and that the 
scores are reasonably reliable and robust. 
Deviations from these assumptions contribute to noise in the system.
Statistical analyses and a list of the critical masses for a variety of academic disciplines 
(not including statistics and operational research) are given in (Kenna and Berche, 2010b).
In the next section, we perform a similar analysis for the statistics and operational research groups submitted to RAE~2008.

%%%%%%%%%%%%%%%%%%%%%%%%%%%%%%%%%%%%%%%%%%%%%%%%%%%%%%%%%%%%%%%%%%%
\section{Statistical analysis of statistics and operational research groups}
%%%%%%%%%%%%%%%%%%%%%%%%%%%%%%%%%%%%%%%%%%%%%%%%%%%%%%%%%%%%%%%%%%%

The  Statistics and Operational Research UOA at RAE~2008 included theoretical, applied and methodological approaches to statistics, probability and  operational research. 
There were 30 submissions comprising 388.8 individuals (with fractions corresponding to part-time staff)
and group sizes ranged from $N=2$ to $N=30$, with mean group size 13.
We find it useful to compare to the Applied Mathematics UOA because of the high degree of overlap between the two disciplines. 
There were 45 submissions in applied mathematics entailing 850.05 individuals
in groups of size $N=1$ to $N=80.3$ with mean group size $18.9$.
The 30 submissions for statistics and operational research are listed in Table~1. Also listed are 
 the numbers of staff submitted and the resultant quality score.

\begin{table}
\caption{Universities which submitted to the Statistics and Operational Research UOA at RAE 2008, listed alphabetically together with the numbers of staff submitted $N$ and quality measurements $s$.}
\centering
\begin{tabular}{|r|l|r|r|} \hline \hline
Index & University &  $N$ & $s$      \\
                                        &       &     &          \\
\hline
 1  & Bath                   	              &	15.00 & 42.14     \\
 2  & Bristol	              	              &	23.00 & 48.57     \\
 3  & Brunel                           & 10.00 & 35.71     \\
 4  & Cambridge	                          & 16.00 & 52.86     \\
 5  & Durham	                           	  &	11.60 & 30.71     \\
 6  & Glasgow	                          	  &	13.00 & 35.71     \\
 7  & Greenwich	                          &	 2.00 & 22.86     \\
 8  & Imperial           	                &	13.90 &	50.00     \\
 9  & Joint submission: Edinburgh 
     \& Heriot-Watt                          &	30.00 &   31.43 \\
10  & Kent	                         &	12.00 &	43.57     \\
11  & Lancaster              	             &	21.65 &	39.29     \\
12  & Leeds	                               &	11.00 &	46.43     \\
13  & Liverpool 	                        & 	 5.00 &	22.14     \\
14  & London Metropolitan                  &	 4.00 &	19.29     \\
15  & London School of Economics 
     \& Political Science  	               &	13.00 &	37.14     \\
16  & Manchester             	             &	10.90 & 39.29     \\
17  & Newcastle              	             &	13.00 &	35.00     \\
18  & Nottingham	                        &	 9.00 &	45.71     \\
19  & Open University                    &	 7.00 &	33.57     \\
20  & Oxford	                           &	24.50 &	62.86     \\
21  & Plymouth              	             &	 4.00 &	19.29     \\
22  & Queen Mary                          &	 8.20 &	29.29     \\
23  & Reading	                             &	 7.70 &	25.71     \\
24  & Salford	                           &	 9.80 &	22.86     \\
25  & Sheffield	                          &	10.70 &	35.71     \\
26  & Southampton	                       &	28.90 &	40.71     \\
27  & St Andrews	                          &	 7.00 & 36.43     \\
28  & Strathclyde	         	         &	10.33 &	29.29     \\
29  & University College London	         &	10.50 &	32.86     \\
30  & Warwick	                           &	24.00 &	48.57     \\
\hline
  \multicolumn{2}{|l|}{ Mean:}                    &  36.50 & 12.96     \\
\hline \hline
\end{tabular}
\end{table}

%...........................................................................
%\begin{figure}[!ht]
\begin{figure}[!t]
\begin{center}
%\hspace{-1cm}
\includegraphics[width=0.45\columnwidth, angle=0]{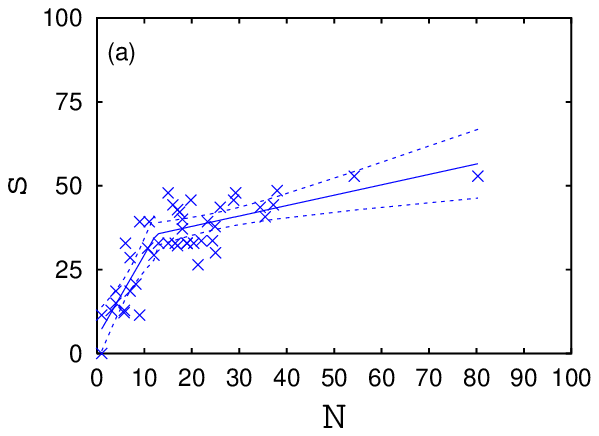}
\includegraphics[width=0.45\columnwidth, angle=0]{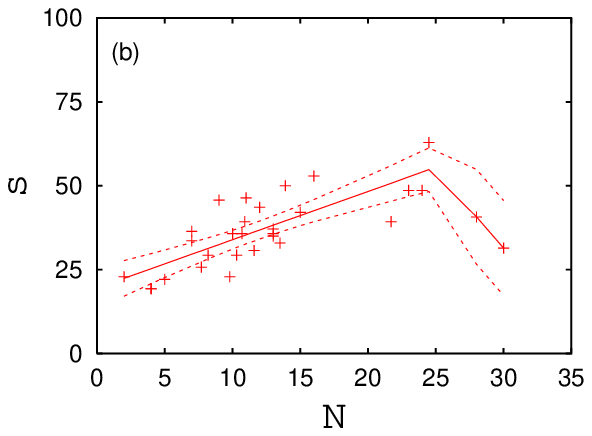}
\caption{Panel (a) depicts quality of research versus quantity of researchers for the Applied Mathematics UOA at RAE~2008 together with the best fit to model (\ref{Nc}) and
95\% confidence interval.
Panel (b) is the equivalent plot for {\emph{all}\/} statistics and operational research groups.}
\end{center}
\end{figure}
%..........................................................................

In Fig.~1(a), we plot RAE-measured quality scores against group quantity for the Applied Mathematics UOA. 
As expected from (\ref{Nc}), research quality indeed tends to increase linearly with group size $N$ up to a 
breakpoint, estimated at  $N_c=12.5 \pm 1.8$ and which splits the 45 research teams into 16 small/medium groups and 29 large ones.
The coefficient of determination is measured to be $R^2=0.74$ and the data passes the Kolmogorov-Smirnov normality test.
The $P$ value for the null hypothesis that there is no underlying correlation between quality and quantity is less than $0.001$, indicating that this can be rejected.
The presence of the breakpoint is evidenced by the $P$ value for the hypothesis that the slopes to the left and right coincide. This is also less than $0.001$, so the
hypothesis can be rejected. The dependency of quality  on quantity continues at a reduced level to the right of the breakpoint as the $P$ value for 
vanishing slope to the right is $0.001$.

In Fig.~1(b), the equivalent full data set for the Statistics and Operational Research UOA is plotted, and the difference between this data set and
that for Applied Mathematics is immediately apparent. 
A correlation between quality and quantity is visible up to about $N=24$, beyond which there are only two data points.
However, the relatively high value of the breakpoint compared to that of applied mathematics (expected to be a closely related
discipline) gives cause for concern, as does the negative slope on the right. No other discipline analysed in (Kenna and Berche, 2010b) 
 exhibited such a phenomenon and
this concern is the reason for the omission of an analysis of statistics and operational research there.

However, closer inspection of the data reveals that the submission with the largest $N$ value, and that corresponding to the rightmost point 
in Fig.~1(b) is in fact a joint submission between Edinburgh and Heriot-Watt universities. This was the only joint submission in this subject area.
Arguing that this submission does not represent a single cohesive ``research group'' in the same spirit as the others in the discipline, we may
consider the corresponding data point to be an outlier and omit it from the analysis.

%...........................................................................
%\begin{figure}[!ht]
\begin{figure}[!t]
\begin{center}
%\hspace{-1cm}
\includegraphics[width=0.45\columnwidth, angle=0]{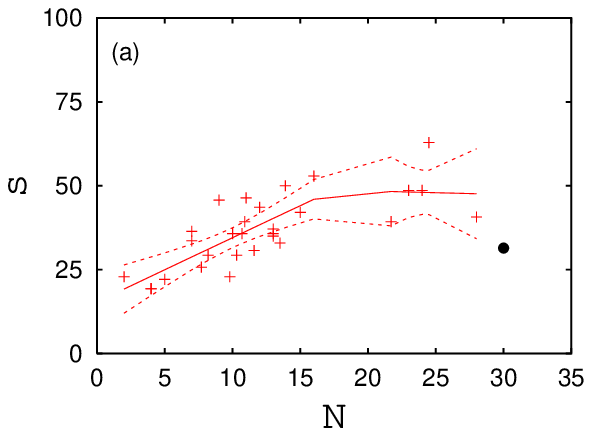}
\includegraphics[width=0.45\columnwidth, angle=0]{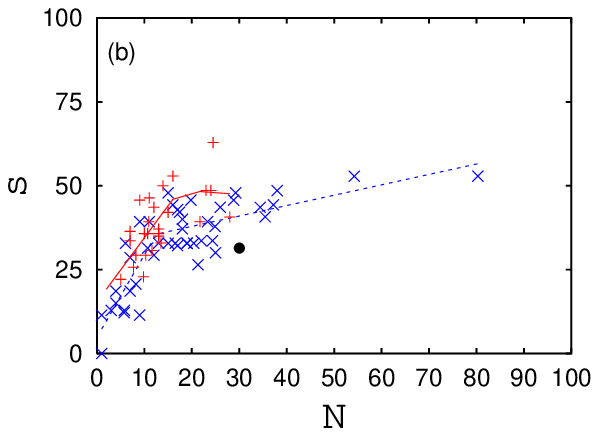}
\caption{(a) The same data as in Fig.~1(b), but omitting that corresponding to the joint submission of Edinburgh and Heriot-Watt universities
(which corresponds to the black disc) from the fitting procedure. 
(b) A comparison between statistics \& operational research (``\textcolor{red}{+}'' symbols and solid  line (red online)) and applied mathematics (``\textcolor{blue}{$\times$}'' symbols and dashed line (blue online)). }
\end{center}
\end{figure}
%..........................................................................

The remaining data are depicted by  crosses (in red online) in the quality versus quantity plot of Fig.~2(a),
in which the Edinburgh/Heriot-Watt  datum is represented by a black circle.
The solid line  is a piecewise linear regression  to the data for which the  dashed curves  represent the 95\%  confidence interval.
One finds a breakpoint at $N_c = 17.4 \pm 5.6$. 
The coefficient of determination is $R^2=0.60$ and the data passes the Kolmogorov-Smirnov normality test.
As for applied mathematics, the $P$ value for the absence of a correlation between quality and quantity is less than $0.001$.
However, unlike applied mathematics,  the P value for the absence of correlation between $s$ and $N$ for large groups  is $0.9$, so this hypothesis cannot be dismissed.
This observation is consistent with the results for other disciplines presented in (Kenna and Berche, 2010b), where we found that research quality tends to saturate in large groups provided $N_k > 7$. 
Also unlike in applied mathematics, the $P$ value for the coincidence of slopes on either side of the transition $N_c$ is  $0.2$ and the corresponding hypothesis cannot be safely disgarded.
We nonetheless arrive at the estimate for the lower critical mass for statistics and operational research given in Eq.(\ref{result}).
This result appears reasonable as it is close to that of applied mathematics, which is $N_k= 6 \pm 1$.

Of course it is possible to fit to other ans{\"{a}}tze, such as polynomials, log-linear curves and power-laws. The results of such fits 
are given in Table~2. Unlike our model (\ref{Nc}) however, these ans{\"{a}}tze are not based on microscopic considerations and interpretation
of, and comparisons between the corresponding results are more difficult. Indeed, we know of no way to extract critical masses from these procedures.
Edinburgh and Heriot-Watt Universities also submitted jointly to the Applied Mathematics UOA at RAE 2008. 
We find that the results of the fit to (\ref{Nc}) are not appreciably affected by removing the datum corresponding to this joint submission.
Notwithstanding this, the statistics reported in Table~2 for applied mathematics correspond to the data set with 
Edinburgh/Heriot-Watt removed. These results are almost identical to those presented in (Kenna and Berche, 2010a;2010b) for the full data set.

\begin{table}
\caption{Results for the model (\ref{Nc}) and for alternative fitting ans{\"{a}}tze.
The Edinburgh/Heriot-Watt joint submissions have been removed from analyses of both disciplines.}
\centering
\begin{tabular}{|l|l|c|c|} \hline \hline
Ansatz for $\langle{s(N)}\rangle$& Parameter      & Applied                         &  Statistics  \&              \\
                                 & and            & mathematics                     &  operational                 \\
                                 &  $R^2$-value   &                                 &  research                    \\
                    \hline
$a_1 + b_1 N$ if $N \le N_c$     & $a_1$          & $5 \pm 4$                       &   $15 \pm	5$                \\
$a_2 + b_2 N$ if $N \ge N_c$     & $b_1$          & $2.5 \pm 0.6 $                  &   $1.9 \pm 	0.5  $          \\
                                 & $a_2$          & $32	\pm 13 $                   &   $51	\pm 35 $              \\
                                 & $b_2$          & $0.4 \pm 0.1  $                             &   $0	\pm 2 $                 \\
                                 & $N_c$          & $12 \pm 2    $                  &   $18	\pm 6 $               \\
                                 & $R^2$          & $74.2 $                         &   60.3\%                     \\
\hline
$A_0 + A_1 N + A_2N^2$           & $A_0$          & $13 \pm 3$                      &   $12 \pm	6$                \\
                                 & $A_1$          & $1.5	0.3 $                      &   $2.9 \pm 	0.9  $          \\
                                 & $A_2$          & $-0.012 \pm 0.003$              &   $-0.059	\pm 0.027 $       \\
                                 & $R^2$          & $67.2\% $                       &   59.9\%                     \\
\hline
$B_0 + B_1 N + B_2 N^2 + B_3N^3$ & $B_0$          & $8 \pm 4$                       &   $17 \pm	9 $               \\
                                 & $B_1$          & $2.4 \pm 0.5$                   &   $1	\pm 3 $                 \\
                                 & $B_2$          & $-0.05 \pm 0.02$                &   $0.10	\pm 0.2 $           \\
                                 & $B_2$          & $0.0003 \pm 0.0002$             &   $-0.004 \pm	0.005 $       \\
                                 & $R^2$          & 70.6\%                          &   61.1\%                     \\
\hline
$C_0 + C_1 N^{C_2}$              & $C_0$          & $-112 \pm 231$                  &   $-15   \pm	75  $           \\
                                 & $C_1$          & $115 \pm 227$                   &   $ 27	  \pm	66 $            \\
                                 & $C_2$          & $0.1 \pm 0.2$                   &   $ 0.3	\pm	0.5 $           \\
                                 & $R^2$          & 72.5\%                          &   57.5\%                     \\
\hline
$D_0 + D_1 \ln{(N+D_2)}$         & $D_0$          & $-4  \pm 10$                    &   $ -16	\pm 44  $           \\
                                 & $D_1$          & $14  \pm 3$                     &   $ 20	\pm	13  $             \\
                                 & $D_2$          & $0.9 \pm 1.5$                   &   $ -4	\pm	8   $             \\
                                 & $R^2$          & 72.8\%                          &   57.9\%                     \\
                    \hline \hline
\end{tabular}
\end{table}

To further compare statistics and operational research to applied mathematics, we plot the sets of data corresponding to both UOA's in Fig.~2(b)
together with the fits coming from the model (\ref{Nc}).
The similarities in their critical masses are evident, as are the similarities between slopes of the piecewise linear fits, although that for 
statistics and operational research is shifted slightly above that for applied mathematics, indicating a consistently better average performance for comparably sized groups or 
problems with the RAE due to the absence of a systematic approach to normalize scores between disciplines.
We believe the latter is the more likely scenario.
In any case, it is clear that in comparison to applied mathematics, there are relatively few statistics and operational research teams in the UK and, of those, there are
even fewer which are supercritical (and therefore operating with sufficient resources) in size.
This suggests that greater investment in this subject area is required to achieve optimal research efficiency.

%...........................................................................
%\begin{figure}[!ht]
\begin{figure}[!t]
\begin{center}
%\hspace{-1cm}
\includegraphics[width=0.45\columnwidth, angle=0]{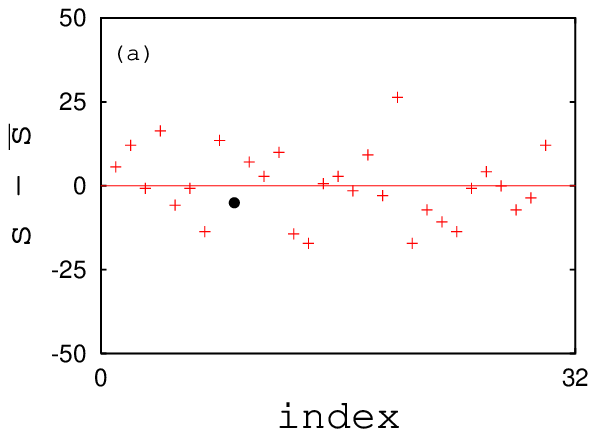}
\includegraphics[width=0.45\columnwidth, angle=0]{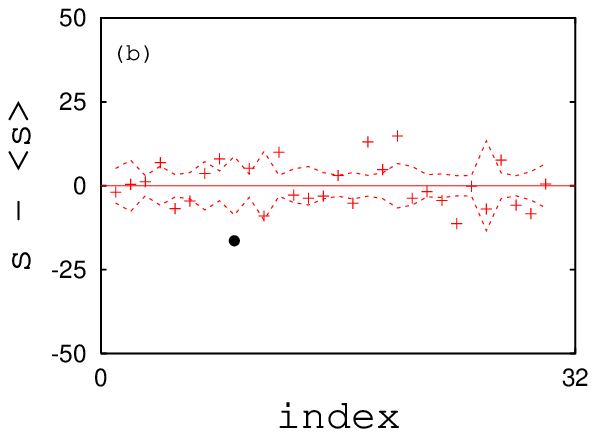}
\caption{(a) Quality measurements normalised to the overall mean for statistics and operational research and (b) renormalised to the expectation values $\langle{s}\rangle$ given in 
Eq.(\ref{Nc}). The tighter distribution of the data about the line in (b) demonstrates the validity of the model. In both plots, the abscissae index the universities listed alphabetically in Table~1.}
\end{center}
\end{figure}
%..........................................................................

To illustrate the superiority of the model over the alternative idea that there is no relationship between quality and quantity in research, we plot 
in Fig.~3 the deviations of the data from the predictions coming from both scenarios.
In each case the data are plotted against the index values listed in Table~1, which correspond to an alphabetical ordering of the institutes which 
submitted to the Statistics and Operational Research UOA. 
In Fig.~3(a), the differences between the quality scores and the mean quality value of the 30 research groups are plotted. 
The range and standard deviation corresponding to this plot are 43.6 and 10.5 respectively (43.6 and 10.7 if Edinburgh/Heriot-Watt is excluded).
In Fig.~3(b), the deviations from the expectation values coming from the model (\ref{Nc}) are plotted. The range and standard deviation
associated with this plot (excluding Edinburgh/Heriot-Watt) are 26.1 and 6.7, respectively.
The tighter distribution of the data in Fig.~3(b) over Fig.~3(a) illustrates the validity of the model.

Plots of the type given in Fig.~3(a) form the basis on which research groups are ranked post RAE, with teams above and below the line deemed to be performing above 
and below average, respectively. However, such rankings do not compare like with like as they fail to take size, and hence resources, into account.
We suggest that Fig.~3(b) forms the basis of a better system as in this plot,  performances are compared to the averages for teams of given sizes.  Fig.~3(b) takes
size into account and gives a better indication of which groups are punching above and below their weights.

%%%%%%%%%%%%%%%%%%%%%%%%%%%%%%%%%%%%%%%%%%%%%%%%%%%%%%%%%%%%%%%%%%%
\section{Conclusions}
%%%%%%%%%%%%%%%%%%%%%%%%%%%%%%%%%%%%%%%%%%%%%%%%%%%%%%%%%%%%%%%%%%%

To summarise, we have applied a mean-field inspired model to examine the relationship between the quality of research teams in statistics and operational research
and the quantity of researchers in those teams. 
Our empirical data is taken from the most recent Research Assessment Exercise in the UK. 
We find that, when an outlying amalgamated group is omitted the dependency of quality upon quantity for this subject area is similar to, and consistent with,  a multitude 
of other disciplines which were reported on in (Kenna and Berche, 2010b). 
The model allows the definition of two critical masses for the discipline.
the research quality of small ($N<N_k$) and medium ($N_k \le N < N_c$) teams is strongly dependent on the number of researchers in the group.
Beyond $N_c$, large teams tend to fragment and research quality is no longer correlated with group size. 
The lower critical mass for statistics and operational research is determined to be $N_k = 9 \pm 3$, and the upper value is about twice that.
These values compare satisfactorily to the equivalent for applied mathematics which has $N_k = 6 \pm 1$. 
To further contextualize these values, we quote from Kenna and Berche (2010b) 
 the results $N_k \le 2$ for pure mathematics (a relatively solitary research discipline)
 and $N_k=20 \pm 4$ for medical sciences (a highly collaborative one).

Notwithstanding the fact that some statisticians and operational researchers were submitted to RAE 2008 as part of 
teams in other disciplines such as   business, economics, engineering and epidemiology,
 about a quarter of statistics/operational research {\emph{groups}\/} submitted to RAE are sub-critical, with $N < N_k = 9$, and therefore vulnerable.
These teams need to strive to attain critical mass. 
Of the 29 teams excluding the Edinburgh/Heriot-Watt combination,
only five (17\%) have size above the upper critical mass of $N_c=18$.
%, four of which are located in the southern half of England.
Therefore the majority of statistics and operational research teams within the UK are under-resourced in terms of staff numbers. 
We suggest that to increase research  efficiency for this discipline investment is needed. 
This conclusion parallels that of Smith and Staetsky (2007) for the teaching of statistics in the UK.

%
%%%%%%%%%%%%%%%%%%%%%%%%%%%%%%%%%%%%%%%%%%%%%%%%%%%%%%%%%%%%%%%%%%%
\vspace{1cm}
\noindent
{\bf{References}} 
%%%%%%%%%%%%%%%%%%%%%%%%%%%%%%%%%%%%%%%%%%%%%%%%%%%%%%%%%%%%%%%%%%%

\vspace{0.2cm}
\noindent
Harrison,~M. (2009) 
Does high quality research require critical mass? 
In {\emph{The Question of R\&D Specialisation: Perspectives and Policy Implications}}
(eds D.~Pontikakis, D.~Kriakou and R.~van~Baval), 
pp~57-59.
European Commission: JRC Technical and Scientific Reports.

\vspace{0.2cm}
\noindent
Evidence, a division of Thompson Reuters (Scientific) Ltd. (2010)
{\emph{The future of the UK University research base}\/} (report for Universities UK).

\vspace{0.2cm}
\noindent
Kenna,~R. and Berche,~B. (2010a) 
The extensive nature of group quality.
{\emph{EPL}\/}, {\bf{90}},  58002.

\vspace{0.2cm}
\noindent
Kenna,~R. and Berche,~B. (2010b) 
Critical mass and the dependency of research quality on group size.
{\emph{Scientometrics}\/} DOI: 10.1007/s11192-010-0282-9.

\vspace{0.2cm}
\noindent
Smith,~T.M.F. and Staetsky, L. (2007)
The teaching of statistics in UK universities. 
{\emph{J.~R.~Statist. Soc.~A}\/}, {\bf{170}}, Part 3, pp.~1-42.

\vspace{1cm}
\noindent
{\bf{Acknowledgements}} 
We are grateful to Houshang Mashhoudy and Neville Hunt for inspiring discussions.

\end{document}